\documentclass[a4,12pt]{article}
\begin{document}
\title{A Theorem of Probability                     
          }
\author{Masumi Nakajima \  \\
        \it Department of Economics \  \\
        \it International University of Kagoshima \   \\
        \it Kagoshima 891-0191, JAPAN  \\
        e-mail: nakajima@eco.iuk.ac.jp}
\maketitle
\begin{abstract}
We prove here the above titled theorem the applications of which will be given elsewhere. \\

Key words ; almost sure convergence.

Mathematics Subject Classification ;  60F99. 
\end{abstract}
Let $(\Omega,{\cal F},P)$ be a probability space and $X_n:=X_n(\omega),\ (n=1,2,3,\cdots)$
 be its random variables with $X_n(\omega) \geq 0$ for $\forall\omega \in \Omega $.
\[
{\rm E}[X]:={\rm E}[X(\omega)]:=\int_{\Omega}X(\omega)dP(\omega)
\]
denotes the expectation value(or mean) of the random variable $X=X(\omega)$.
We put $\{ K_m\}_{m=1}^\infty$ to be a natural number sequence with $K_m <K_{m+1}$ and 
$K_m \to +\infty \ ({\rm as} \ m \to +\infty)$.\\
Then we have
\newtheorem{theo}{Theorem}
%
\begin{theo}
If
\begin{eqnarray}
{\rm E}[\sum_{m=1}^\infty X_{K_m + l_m}(\omega)]<+\infty \ for \ 
0<\forall l_m \leq K_{m+1}-K_m , 
\end{eqnarray}
then we have 
\begin{eqnarray}
X_n(\omega) \to 0 \ almost \ surely \ on \ \Omega \ (as \ n \to \infty).  
\end{eqnarray}
\end{theo}

\newtheorem{lem}{Lemma}

{\bf proof} \ The conclusion of the theorem is equivalent to
\begin{eqnarray}
P\{ \omega \in \Omega \ | \ \forall \epsilon >0,\ \exists N \equiv N(\epsilon,\omega),\ 
X_n(\omega)\leq \epsilon \ for \ \forall n \geq N \}=1
\end{eqnarray} 
Firstly we assume our denying of this conclusion, that is, we assume
\begin{eqnarray}
P\{ \omega \in \Omega \ | \ \exists \epsilon(\omega) >0,\ \forall N \in {\bf N},\ 
\exists n \equiv n(\omega) \geq N, \  
X_n(\omega) > \epsilon(\omega) \ \}>0 
\end{eqnarray}
which will lead to a contradiction later.\\
(4) is equivalent to $P\{ A \}>0$ with 
\begin{eqnarray}
A:=\{ \omega \in \Omega \ | \ \exists \epsilon(\omega) >0, 
\exists \{n_j \equiv n_j(\omega) \}_{j=1}^\infty,  n_j \to \infty ,  
X_{n_j}(\omega) > \epsilon(\omega)  \}
\end{eqnarray}
We put the followings:
\begin{eqnarray}
B_n&:=&\{ \omega \in \Omega \ |\ X_n(\omega)>\epsilon(\omega) \} \\
\{m_j\}_{j=1}^\infty&:=&\{ n \in {\bf N} \ |\ P\{B_n\}>0 \} \ \ with \ \ m_j<m_{j+1} \\
C_N&:=&\bigcup_{m_j\geq N}B_{m_j}
\end{eqnarray}
Then we have trivially
\begin{eqnarray}
C_N \supset C_{N+1}
\end{eqnarray}
which shows the existence of 
\begin{eqnarray}
C:=\lim_{N \to \infty}C_N={\limsup}_{j \to \infty}B_{m_j}
\end{eqnarray}
and
\begin{eqnarray}
C=A \ \ \ except \  P-measure \  zero \  set
\end{eqnarray}
Next we divide $\{m_j\}_{j=1}^\infty$ into $\{m_j(k)\}_{j=1}^\infty$'s as follows: 
\begin{eqnarray}
\{m_j\}_{j=1}^\infty=\bigcup_{k \in \Lambda} \{m_j(k)\}_{j=1}^\infty
\end{eqnarray}
with
\begin{eqnarray}
m_j(k)&<&m_{j+1}(k) \\
\{m_j(k)\}_{j=1}^\infty \bigcap \{m_j(l)\}_{j=1}^\infty &=& \emptyset (k \neq l) \\
\sharp \{ \{m_j(k)\}_{j=1}^\infty \bigcap (K_m,K_{m+1}] \} &\leq& 1 \ for \ \forall m \in 
{\bf N} \\
\sharp \Lambda &=& \sharp {\bf N} ( i.e. countablly \ many )
\end{eqnarray}
where $\sharp A $ denotes the number of elements of the set $A$.\\
We also put
\begin{eqnarray}
D_N^{(k)}:=\bigcup_{m_j(k) \geq N}B_{m_j(k)}.
\end{eqnarray}
Then we also have trivially
\begin{eqnarray}
D_N^{(k)} \supset D_{N+1}^{(k)}
\end{eqnarray}
which leads to the existence of
\begin{eqnarray}
D^{(k)}=\lim_{N \to \infty}D_N^{(k)}={\limsup}_{j \to \infty}B_{m_j(k)}
\end{eqnarray}
and
\begin{eqnarray}
C_N=\bigcup_{k \in \Lambda}D_N^{(k)}.
\end{eqnarray}
From (20), it follows that when $N$ tends to $\infty$
\begin{eqnarray}
A=C=\bigcup_{k \in \Lambda}D^{(k)} \ \ \ except \  P-measure \  zero \  set
\end{eqnarray}
which leads to
\begin{eqnarray}
\exists l \in {\bf N} \ \ such \ \ that \ \ P\{D^{(l)}\}>0
\end{eqnarray}
because of $P\{A\}>0$ and $\sharp\Lambda=\sharp{\bf N}.$ \\
We put
\begin{eqnarray}
K:=\int_{D^{(l)}} \epsilon (\omega) dP(\omega).
\end{eqnarray}
Because of the assumption of the theorem (1), that is,
\begin{eqnarray}
{\rm E}[\sum_{j=1}^\infty X_{m_j(l)}(\omega)]<+\infty 
\end{eqnarray}
there exists a natural number $N$ such that
\begin{eqnarray}
{K\over 2}&>&{\rm E}[\sum_{m_j(l)\geq N} X_{m_j(l)}(\omega)] \\
&=& \sum_{m_j(l)\geq N}\int_\Omega X_{m_j(l)}(\omega)dP(\omega) \\
&>& \sum_{m_j(l)\geq N}\int_{B_{m_j(l)}} X_{m_j(l)}(\omega)dP(\omega) \\
&>& \sum_{m_j(l)\geq N}\int_{B_{m_j(l)}} \epsilon(\omega)dP(\omega) \\
&\geq& \int_{D_N^{(l)}} \epsilon(\omega)dP(\omega)
\end{eqnarray}
due to $D_N^{(l)}=\bigcup_{m_j(l) \geq N} B_{m_j(l)}$.
Because of
\begin{eqnarray}
\int_{D_N^{(l)}} \epsilon(\omega)dP(\omega)\geq \int_{D_{N+1}^{(l)}} \epsilon(\omega)dP(\omega)
\to \int_{D^{(l)}} \epsilon(\omega)dP(\omega)=K,
\end{eqnarray}
we have 
\begin{eqnarray}
\int_{D_N^{(l)}} \epsilon(\omega)dP(\omega)
\geq K
\end{eqnarray} 
with sufficiently large $N$.
From (25)$,\cdots,$(31), we have a contradiction:
\begin{eqnarray}
{K \over 2} \geq K .
\end{eqnarray} 
Therefore we cannot have (4) or (5) which means the conclusion of the theorem:(2) or (3).
This completes the proof. 

\end{document}